\documentclass[showpacs, amsmath,amssymb, aps]{revtex4}

\usepackage{tabularx}
\usepackage{graphicx}
\usepackage{float}
\usepackage{dcolumn}
\usepackage{bm}
\usepackage{mathrsfs}
\usepackage{amsmath}
\newtheorem{theorem}{Theorem}[section]
 \newtheorem{definition}[theorem]{Definition}
\newtheorem{lemma}[theorem]{Lemma}

\newtheorem{corollary}[theorem]{Corollary}

\begin{document}


\title{ Explicit Superconic Curves}
\author{ Sunggoo Cho}
\email{sgcho@semyung.ac.kr, oldrock9@hanmail.net}
\affiliation{\small School of Computer Science,
Semyung University, Chechon, Chungbuk 390 - 711, Korea}

\date{\today}

\begin{abstract}

Conics and Cartesian ovals are very important curves in various fields of science.
Also aspheric curves based on conics  are useful in optics.
Superconic curves recently suggested by A. Greynolds  are extensions
of  both conics and Cartesian ovals  and have been applied to
optics while they are not extensions of aspheric
curves based on conics.

In this work, we investigate another kind of superconic curves
that are extensions of not only conics and Cartesian ovals
but also aspheric curves based on conics. Moreover, the superconic
curves are represented in explicit form
while  Greynolds' superconic curves are in implicit form.

\end{abstract}

\pacs{02.30.Gp, 02.60.Cb, 42.15.Dp }

\keywords{Cartesian ovals, aspheric curves, Superconic curves, Optical design,}

\maketitle

\section{INTRODUCTION}

Conics and Cartesian ovals are very useful curves in science\cite{Cun1989, Law1972, Loc1967, Kra2004}.
Especially, in optical design, the conic curves
with curvature $c_0$ and conic constant $K$ are  described in implicit and explicit
form as follows\cite{Sha1997, Gal2000}:
\begin{eqnarray}
(1+K)c_0z^2-2z+c_0y^2&=&0,\label{conicIm}
\end{eqnarray}
and
\begin{eqnarray}
z&=&c_0y^2/\{1+[1-(1+K)c_0^2y^2]^{\frac{1}{2}}\},\label{conicEx}
\end{eqnarray}
usually in $(z, y)$-plane for representing  rotationally symmetric surfaces
about the $z$-axis.  We note that the explicit form like
Eq.(\ref{conicEx}) is essential in optical design since it can be
interpolated even in the region where  the curve  is not defined. For an
example, there is no solution $z$ of a circle of radius 1 when
$y>1$, which is the case when the square root term $1-(1+K)c_0^2y^2$
in the denominator becomes negative. In this case, we may
interpolate the circle as parabola $z=c_0y^2$ by putting the square
root term to be zero for $y>1$. Such interpolations enable us to construct so called
{\em an aspheric curve based on conic}  that  is usually
represented for any $y$ in explicit form \cite{Gal2000} as
\begin{eqnarray}
z=c_0y^2/\{1+[1-(1+K)c_0^2y^2]^{\frac{1}{2}}\}+ \sum_{n=2}^{N} f_{2n}
y^{2n}\label{asp}
\end{eqnarray}
for some positive integer $N$ and coefficients $f_{2n}$'s.

Another useful curves are  Cartesian ovals\cite{Law1972, Loc1967}. A Cartesian oval is defined as
the set of points such that the sum of whose weighted distances from two fixed
points is a constant. In general, Cartesian ovals are quartic equations and they are famous for
their perfectly focusing refraction property in optics\cite{Bev2006, Mic2011,Rud1998, WinMinBeb2005}.

On the other hand, superconic curves\cite{Ala2002} were suggested by A. Greynolds
 as extensions of  conics and Cartesian ovals, which are
expressed in implicit form as
\begin{eqnarray}
c_0Kz^2-2(1+
b_1s^2+b_2s^4+\cdots)z+(c_0+c_1s^2+c_2s^4+\cdots)s^2&=&0,\label{Reynolds}
\end{eqnarray}
for some constants $b_1, b_2, \cdots$ and $c_1, c_2, \cdots$. Here the parameter $s$ is defined by $s^2\equiv z^2+y^2$. Hence
it is obvious that superconic curves are conics if $b_i=0$ and $c_i=0$
for $i=1, 2, \cdots$. If $ b_1\ne 0 $ and $c_1\ne 0$
 while $b_i=0$ and $c_i=0$ for $i=2, \cdots$, they are Cartesian ovals. However, it is
obvious that the superconic curves cannot be  extensions of
 aspheric curves based on conics  of  Eq.(\ref{asp}).  Furthermore,
it seems to be hard to express Eq.(\ref{Reynolds}) in explicit form
since $s^2=z^2+y^2$. In fact, Greynolds seemed to give up a closed-form
explicit representation corresponding to the implicit form as he stated in his
work\cite{Ala2002}.

In this work, we shall investigate another kind of superconic curves,
which  are extensions of both conic and  Cartesian oval curves like Greynolds'.
On the contrary to Greynolds', however,
our superconic curves are not only extensions of aspheric surfaces based on conics of Eq.(\ref{asp})
but also expressed in explicit form.

In this work we are interested in a solution that  passes
through the origin in $(z, y)$-plane among four solutions of the
quartic equation described by a Cartesian oval.
\begin{definition}
An optical solution is defined to be the solution that passes through the origin
among solutions of a Cartesian oval.
\end{definition}

The main strategy of this work is to find the optical solution
 such that it is not only expressed in
an appropriate explicit form whose limit is conic but also
interpolated in the region where the Cartesian oval is not
defined.

Of course, any quartic equations can be solved
explicitly in general by the method of L. Ferrari who is
attributed with the discovery of the explicit solutions to the
quartic equations in 1540, and they are still studied for their diverse solving methods and properties\cite{AbrSte1972, Nic2009, Shm2011}.
However, the usually known general explicit forms for the solutions
of  quartic equations do not seem to be
appropriate for our purpose.

In section II,  we shall decompose a Cartesian oval into a product of two specific quadratic forms,
two solutions of which are the candidates for the optical solution with an appropriate form
for the limit and interpolation.
In section III, we shall investigate the initial criteria for the choice of the
optical solution between the candidates. Moreover,
the continuity and interpolation of the optical solution shall be discussed.
In section IV, we shall show that conics are the limiting cases of  optical solutions
from a different point of view than are usually known in the literature\cite{Loc1967,Vil1, Eve1965}.
    The limiting process in
this section gives us more insights on the relationship between
 optical solutions   and
 conics. The main result of this work about superconic curves comes from the limiting relationship. Finally, in
section V, we shall discuss about a family of curves including  optical solutions and conics
and demonstrate an example.  A lot of computational work are required in this paper
and done using Mathematica.

\section{The candidates for the optical solutions from Cartesian ovals }
\subsection{Motivations and notations for Cartesian ovals}

There are several forms for Cartesian ovals in the literature\cite{Law1972, Loc1967, Bev2006}.
The purpose of this section is to find an appropriate form for the quartic equation of a Cartesian
oval such that not only it yields the optical solution
 but also  the optical solution can be
interpolated in the region where the Cartesian oval is not
defined. Moreover, the optical solution should become  a conic curve
as its special limiting case.

For this purpose, let us describe a Cartesian oval from a physical point of view
that gives us the natural motivation for the Cartesian oval.
Among others, it may be
better to start from the Snell's law of refraction
 on a curve that passes through the origin in the $(z, y)$-plane:
\begin{eqnarray}
n\sin(\phi-\theta) &=& n' \sin (\phi-\theta'),
\end{eqnarray}
where $n$ and $n'$ are the refractive indices, $\theta$ and
$\theta'$ are the angles which the incident  and refracted
ray make with the $+z$ axis respectively, and $\phi$ is defined by $\tan \phi
\equiv -dz/dy$.
Now let $(z_i, y_i)$ be the position of a point light source
and $(z_o, y_o)$ be a focusing point.
 Also define  signs $s_i,
s_o$ as follows: if $z_i<z$, $s_i = -1$, otherwise $s_i=+1$.
Similarly, if $z_o<z$, $s_o = -1$,  otherwise $s_o=+1$.

Then the law of refraction on a perfectly focusing curve\cite{Vil1} may be
written as
\begin{eqnarray}
\frac{dz}{dy} &=& (-1) \frac{n' \sin \theta' - n \sin
\theta}{n'\cos \theta' - n \cos \theta}\\
&=& (-1)\frac{s_on'(y_o - y)\sqrt{(z_i - z)^2+(y_i - y)^2} -
s_in(y_i - y)\sqrt{(z_o - z)^2+(y_o - y)^2}}{s_on'(z_o -
z)\sqrt{(z_i - z)^2+(y_i - y)^2} - s_in(z_i - z)\sqrt{(z_o -
z)^2+(y_o - y)^2}},\nonumber
\end{eqnarray}
which is the first order exact differential equation
whose solution is given by
\begin{eqnarray}
 s_in\sqrt{(z-z_i)^2+(y-y_i)^2}-s_on'\sqrt{(z-z_o)^2+(y-y_o)^2}&=&n\kappa,
 \label{eq8Sol}
\end{eqnarray}
where $\kappa$ is a constant to  be determined by the initial
position of the curve.

Eq.(\ref{eq8Sol}) is one of many different forms of
Cartesian ovals.
Without loss of generality, we assume that $y_i = y_o = 0$ and
define $m\equiv n'/n >0$. Thus a Cartesian oval is described as the
set of points $(z,y)$
 satisfying the following equation:
\begin{eqnarray}
s_i[(z - z_i)^2+ y^2]^{\frac{1}{2}}-ms_o[(z - z_o)^2+
y^2]^{\frac{1}{2}}&=& \kappa.\label{eq1}
\end{eqnarray}
We also assume that  it
passes through the origin $(0,0)$. Then $\kappa=z_i-mz_o$.

Now let us define  $ \eta_i \equiv \epsilon /z_i$ and $ \eta_o
\equiv \epsilon /z_o$ for a constant $\epsilon>0$ since the
notations are convenient for later use. Then Eq.(\ref{eq1}) may be
written as
\begin{eqnarray}
\eta_o[(\eta_i z-\epsilon )^2+ (\eta_i y)^2]^{\frac{1}{2}}
 -m \eta_i[(\eta_o z-\epsilon)^2+ (\eta_o y)^2]^{\frac{1}{2}}&=&\eta_i\eta_o \kappa \equiv k\label{eq2}
\end{eqnarray}
with $k=\epsilon(\eta_o-m\eta_i)$.

In the next subsection, we shall decompose a Cartesian oval described by Eq.(\ref{eq2}) into
the product of two quadratic factors, two solutions of which are the candidates for the optical solution.

\subsection{The candidates for the optical solution}

Now we are interested in a quartic equation described by Eq.(\ref{eq2}).
Thus in order to avoid those cases for which  Eq.(\ref{eq2}) becomes quadratic,
we assume that  $k \ne 0$, $ m \ne 1$, $\eta_i \ne \eta_o$, $\eta_i\ne0$ and $\eta_o\ne 0$.

Now we define $ x \equiv
\eta_i \eta_o z$. Then  Eq.(\ref{eq2}) may be rewritten as follows.
\begin{eqnarray}
s_o[(x-\eta_o\epsilon)^2+ (\eta_i\eta_o y)^2]^{\frac{1}{2}} &=&k +m
s_i[(x-\eta_i\epsilon)^2+ (\eta_i\eta_o
y)^2]^{\frac{1}{2}}.\label{quartic}
\end{eqnarray}

Now it is easy to see that Eq.(\ref{quartic}) can be transformed into the
following specific form  by
squaring twice to remove the two root terms
\begin{eqnarray}
(b_2x^2+b_1x+b_0)^2&=&a_2x^2+a_1x+a_0 \label{eq10}.
\end{eqnarray}
where
\begin{eqnarray}
a_2&\equiv& 4k^2m^2,\\
a_1&\equiv&-2a_2\epsilon\eta_i,\nonumber\\
a_0&\equiv&a_2 \epsilon^2\eta_i^2+a_2 \eta_i^2\eta_o^2y^2,\nonumber
\end{eqnarray}
and
\begin{eqnarray}
b_2&\equiv&1-m^2,\\
b_1&\equiv&-2\epsilon(\eta_o-m^2\eta_i),
\nonumber\\
b_0 &\equiv&
\epsilon^2(\eta_o^2-m^2\eta_i^2)-k^2+(1-m^2)\eta_i^2\eta_o^2y^2.\nonumber
\end{eqnarray}
Then  we add both sides with
$(b_2x^2+b_1x+b_0+\lambda)^2-(b_2x^2+b_1x+b_0)^2$ to obtain
\begin{eqnarray}
(b_2x^2+b_1x+b_0+\lambda)^2 &=&(a_2+2\lambda
b_2)(x+\frac{a_1+2\lambda b_1}{2(a_2+2\lambda b_2)})^2,\label{eq11}
\end{eqnarray}
where $\lambda$ is supposed to satisfy, unless $a_2+2\lambda b_2=0$,
\begin{eqnarray}
(a_0+2\lambda b_0)+\lambda^2&=&\frac{(a_1+2\lambda
b_1)^2}{4(a_2+2\lambda b_2)}.
\end{eqnarray}
Thus we have a resolvent cubic equation\cite{AbrSte1972} in $\lambda$
\begin{eqnarray}
a\lambda^3+b\lambda^2+c\lambda+d=0,\label{lambda}
\end{eqnarray}
where the coefficients, after divided by $8$, are given by
\begin{eqnarray}
a &\equiv& b_2 ,\label{abcd}\\
b &\equiv& (a_2+4b_2b_0-b_1^2)/2,\nonumber\\
c &\equiv& a_2b_0+a_0b_2-a_1b_1/2,\nonumber\\
d &\equiv& a_2a_0/2-a_1^2/8.\nonumber
\end{eqnarray}

In terms of $m, \eta_i$ and $\eta_o$, Eq.(\ref{abcd}) is written as
\begin{eqnarray}
a&=& 1-m^2 ,\label{coeff}\\
b &=&2 \epsilon^2(m-1)[\eta_o^2(1+m)+2\eta_i^2m^2(1+m)-2\eta_i\eta_om(1+2m)]+2\eta_i^2\eta_o^2(m^2-1)^2y^2,\nonumber\\
c &=& 4\epsilon^4m^2(1-m)\eta_i(\eta_i-2\eta_o+m\eta_i)(\eta_o-m\eta_i)^2+8\epsilon^2m^2(1-m^2)\eta_i^2\eta_o^2(\eta_o-m\eta_i)^2y^2,\nonumber\\
d &=& 8\epsilon^4m^4\eta_i^2\eta_o^2(\eta_o-m\eta_i)^4y^2.\nonumber
\end{eqnarray}

Using the root $\lambda$ of Eq.(\ref{lambda}),  Eq.(\ref{eq11}) becomes
\begin{eqnarray}
0&=& \{b_2x^2+[b_1+(a_2+2\lambda
b_2)^{\frac{1}{2}}]x+b_0+\lambda+\frac{a_1+2\lambda
b_1}{2(a_2+2\lambda
b_2)^{\frac{1}{2}}}\}\label{eq18}\\
&&\times \{b_2x^2+[b_1-(a_2+2\lambda
b_2)^{\frac{1}{2}}]x+b_0+\lambda-\frac{a_1+2\lambda
b_1}{2(a_2+2\lambda
b_2)^{\frac{1}{2}}}\}\nonumber\\
&\equiv&(b_2x^2+p_{+}x+q_{+})(b_2x^2+p_{-}x+q_{-}),\nonumber
\end{eqnarray}
where
\begin{eqnarray}
p_{\pm}&\equiv& b_1\pm (a_2 + 2\lambda b_2)^{\frac{1}{2}},\\\label{pqpm}
q_{\pm}&\equiv&b_0+\lambda \pm \frac{a_1+2\lambda b_1}{2(a_2+2\lambda b_2)^{\frac{1}{2}}}.\nonumber
\end{eqnarray}

If we suppose that $p_1, p_2 \ne 0$ and replace $x$ by $z$ in
Eq.(\ref{eq18}), we have the generic form for the quartic equation
\begin{eqnarray}
(A_{+}z^2-2z+B_{+})(A_{-}z^2-2z+B_{-})=0,\label{product}
\end{eqnarray}
where
\begin{eqnarray}
A_{\pm}&\equiv&\frac{-2b_2\eta_i\eta_o}{p_{\pm}}\,\,\,\mbox{and}\,\,\,B_{\pm}\equiv\frac{-2q_{\pm}}{\eta_i\eta_op_{\pm}}. \label{eqA}\\
\nonumber
\end{eqnarray}

Here the  two solutions
$z=B_{\pm}/[1-(1-A_{\pm}B_{\pm})^{\frac{1}{2}}]=[1+(1-A_{\pm}B_{\pm})^{\frac{1}{2}}]/A_{\pm}$
 are excluded since they do not pass
 through the origin, which may be easily seen from the fact that $A_{\pm}$  are finite under the
 condition that $p_{\pm} \ne 0$.

Thus, we  obtain the two candidates for the optical solution from the quartic
equation in Eq.(\ref{product}).

\begin{lemma}
Let $p_{\pm} \ne 0$, $a_2+2\lambda b_2 \ne 0$ and $(1-A_{\pm}B_{\pm})^{\frac{1}{2}} \ge 0$. Then the optical solution of a Cartesian oval in Eq.(\ref{eq2}) is given by
one of the followings:
\begin{subequations}
\label{eq:one}
\begin{eqnarray}
z&=&\frac{B_{+}}{1+(1-A_{+}B_{+})^{\frac{1}{2}}},\label{SolB}\\
z&=&\frac{B_{-}}{1+(1-A_{-}B_{-})^{\frac{1}{2}}}.\label{SolD}
\end{eqnarray}
\end{subequations}
\end{lemma}

In the next section, we shall investigate the criteria to choose the optical solution between Eq.(\ref{SolB}) and Eq.(\ref{SolD}).

\section{The optical solutions}

\subsection{The criteria for the optical solutions }

In order to calculate the optical solution,
we  need to determine $\lambda$ satisfying the cubic equation in Eq.(\ref{lambda}).
If we introduce parameters
\begin{eqnarray}
q &\equiv& (3ac-b^2)/(3a^2),\label{eqdelta}\\
r&\equiv& (9abc-27a^2d-2b^3)/(27a^3),\nonumber\\
\triangle &\equiv& r^2+\frac{4}{27}q^3,\nonumber
\end{eqnarray}
the cubic equation Eq.(\ref{lambda}) can be solved easily by the
well-known method. That is, if $\triangle
> 0$, there are two complex and one real roots, where the real root
is expressed as
\begin{eqnarray}
\lambda
&=&-\frac{b}{3a}+(\frac{r+\sqrt{\triangle}}{2})^\frac{1}{3}+(\frac{r-\sqrt{\triangle}}{2})^\frac{1}{3}.
\label{eq16}
\end{eqnarray}
If $\triangle < 0$, define $\rho \equiv
(-\frac{q^3}{27})^{\frac{1}{2}}$ and $\cos \theta \equiv r/(2\rho)$.
Then there are three real roots given by
\begin{subequations}
\label{eq20}
\begin{eqnarray}
\lambda &= & -\frac{b}{3a}+2
\rho^{\frac{1}{3}}\cos(\frac{\theta}{3}),\label{eq20a}\\
\lambda &= & -\frac{b}{3a}+2 \rho^{\frac{1}{3}}\cos(\frac{\theta-
2\pi}{3}),\label{eq20b}\\
 \lambda &= & -\frac{b}{3a}+2
\rho^{\frac{1}{3}}\cos(\frac{\theta+ 2\pi}{3}).\label{eq20c}
\end{eqnarray}
\end{subequations}


We note that $\triangle $ is a sextic polynomial in $y$ as follows.
\begin{eqnarray}
\triangle  &=& D_6 y^6 + D_4 y^4 +D_2 y^2 + D_0,\label{triangle}
\end{eqnarray}
where
\begin{eqnarray}
D_6 &\equiv& -\frac{256}{27} \epsilon^6\eta_i^6(\eta_i - \eta_o)^2
\eta_o^6m^4(\eta_o - \eta_i m)^4,\label{triangle_coeff}\\
D_4 &\equiv& \frac{64 \epsilon^8 \eta_i^4 \eta_o^4 m^4(\eta_o-
\eta_i m)^4 D_4^*}{27(1 - m)^2(1 + m)^2},\nonumber\\
D_2 &\equiv& \frac{128 \epsilon^{10} \eta_i^2 \eta_o^2 m^4(\eta_o-
\eta_i m)^4 D_2^*}{27(1 - m)(1 + m)^4},\nonumber\\
D_0&\equiv&-\frac{64}{27}\frac{\epsilon^{8}m^4k^4}{(1+m)^4}[(1+m)\eta_i^2- 2\eta_i\eta_o]^2[(1+m)\eta_o^2-2m\eta_i\eta_o]^2,\nonumber
\end{eqnarray}
and
\begin{eqnarray}
D_4^*  &=&\eta_o^4(8 + 20 m^2 - m^4) +
    4\eta_i\eta_o^3(-4 - 4 m - 15 m^2 - 5 m^3 + m^4)\label{triangle_coeff_star}\\
   && +    2\eta_i^2\eta_o^2(2 + 20 m + 37 m^2 + 20 m^3 + 2 m^4) -
    4\eta_i^3\eta_o(-1 + 5 m + 15 m^2 + 4 m^3 + 4 m^4) \nonumber\\
    &&+\eta_i^4(-1 + 20 m^2 + 8 m^4),\nonumber\\
D_2^*  &=& -2\eta_o^6(1 + m)^3 + 2\eta_i^6 m^2(1 + m)^3 +
    2\eta_i\eta_o^5(1 + m)^2(2 + 6 m + m^2) \\
    &&-
    2\eta_i^5\eta_o m(1 + m)^2(1 + 6 m + 2 m^2) +
    \eta_i^4\eta_o^2(1 + 11 m + 29 m^2 + 39 m^3 + 18 m^4 - 2 m^5) \nonumber\\
    &&+
    4 \eta_i^3\eta_o^3(-1 - 2 m + m^2 - m^3 + 2 m^4 + m^5)\nonumber\\
    && -
    \eta_i^2\eta_o^4(-2 + 18 m + 39 m^2 + 29 m^3 + 11 m^4 + m^5).\nonumber
\end{eqnarray}
It follows then  that
\begin{eqnarray}
\triangle&=&-\frac{64}{27}\frac{\epsilon^{8}m^4k^4}{(1+m)^4}\sigma^2\delta^2 <0,\nonumber
\end{eqnarray}
at $y=0$ unless $\sigma \cdot\delta=0$, where $\sigma$ and $\delta$ are defined as
\begin{eqnarray}
\sigma &\equiv&  (1+m)\eta_i^2- 2\eta_i\eta_o,\\
\delta &\equiv&  (1+m)\eta_o^2-2m\eta_i\eta_o. \nonumber
\end{eqnarray}

Thus unless $\sigma \cdot\delta=0$, any real root $\lambda$ in Eq.(\ref{eq20})
may be used  to calculate the solutions of a Cartesian oval.

We shall show  first that
the choice of the optical solution between Eq.(\ref{SolB}) and Eq.(\ref{SolD}) is determined by the initial choice of the root $\lambda$  of Eq.(\ref{lambda}).

Now, under the assumption that $\sigma\ne 0$ and $\delta \ne 0$,
let us find the roots of Eq.(\ref{lambda}) when $y \rightarrow 0$.
In fact, from Eq.(\ref{coeff}), it follows that $d\rightarrow 0$ as $y\rightarrow
0$. Thus the cubic equation of Eq.(\ref{lambda}) in $\lambda$ becomes $
\lambda(a\lambda^2+b\lambda+c)=0$. Hence if we expand $
\lambda=\lambda_0 + \lambda_1 y +\lambda_2 y^2 + \cdots $ when $y$ is small,
it is not difficult to find the three real roots
$\lambda_0$.
\begin{subequations}
\label{eq27}
\begin{eqnarray}
\lambda_{01} &=& 0, \label{lambda01}\\
\lambda_{02}&=& 2k^2,\label{lambda02}\\
\lambda_{03}&=&
\frac{2\epsilon^2m^2\eta_i(\eta_i-2\eta_o+m\eta_i)}{1+m}=\frac{2\epsilon^2m^2\sigma}{1+m}.\label{lambda03}
\end{eqnarray}
\end{subequations}
Also, it is trivial to see that $\lambda_1 = 0$ for each case.
Moreover,
 we may calculate $\lambda_2$ corresponding to
Eqs.(\ref{lambda01})$\sim$(\ref{lambda03}) respectively as follows.
\begin{subequations}
\label{eq28}
\begin{eqnarray}
\lambda_{21}&=&\frac{2m^2\eta_i\eta_o^2(\eta_o-m\eta_i)^2}
{(m-1)(\eta_i-2\eta_o+m\eta_i)},\label{lambda21}\\
\lambda_{22}&=&\frac{2\eta_i^2\eta_o(\eta_o-m\eta_i)^2}{(m-1)(\eta_o-2m\eta_i+m\eta_o)},\label{lambda22}\\
\lambda_{23}&=&\frac{-2\eta_i\eta_o(\eta_i-\eta_o)^4m^2(1+m)}
{(m-1)(\eta_i-2\eta_o+m\eta_i)(\eta_o-2m\eta_i+m\eta_o)},\label{lambda23}
\end{eqnarray}
\end{subequations}
where we note that the denominators do not vanish since $\sigma\ne 0$ and $\delta \ne 0$.
\begin{lemma}
Let us assume that $\sigma\ne 0$ and $\delta \ne 0$.

1. If (1)
 $\lambda=\lambda_{01}$ or $\lambda=\lambda_{02}$ at $y=0$ and $k>0$, or
 (2) $\lambda=\lambda_{03}$  at $y=0$ and $\eta_i >\eta_o$, the optical solution is Eq.(\ref{SolB}).

2. If (1) $\lambda=\lambda_{01}$ or $\lambda=\lambda_{02}$ at $y=0$ and $k<0$, or
 (2) $\lambda=\lambda_{03}$  at $y=0$ and $\eta_i <\eta_o$, the optical solution is Eq.(\ref{SolD}).
\end{lemma}
{\bf Proof:} For small $y$, we put   $B\approx B_0 + B_2 y^2 $ for
$B=B_{\pm}$.
Let us consider the case when $\lambda_0 =\lambda_{01}$ or $\lambda_0=\lambda_{02}$.  Then if $k>0$,
a little bit lengthy but straightforward calculation
using $\lambda_{21}$ and $\lambda_{22}$ in  Eq.(\ref{eq28}) respectively shows that for small $y$
\begin{eqnarray}
B_{+}&\approx&\frac{ \eta_i-m\eta_o}{ \epsilon
(1-m)}y^2,\label{approxiB}
\end{eqnarray}
 while the constant term of $B_{-}$ in $y$ does not vanish. Thus the optical solution is
Eq.(\ref{SolB}) since  it passes through the origin.  If $k<0$, however, we have
\begin{eqnarray}
B_{-}&\approx&\frac{ \eta_i-m\eta_o}{ \epsilon
(1-m)}y^2,\label{approxiD}
\end{eqnarray}
while the constant term of  $B_{+}$ in $y$ does not vanish. Thus the
 optical solution is Eq.(\ref{SolD}).
On the other hand, let us suppose that $\lambda_0=\lambda_{03}$. Then if $\eta_i>\eta_o$, we have the same result
as that of Eq.(\ref{approxiB}) for $B_{+}$ by a straightforward calculation.
Thus the optical solution is Eq.(\ref{SolB}). If $\eta_i<\eta_o$, however, we have
the same result as that of Eq.(\ref{approxiD}) for $B_{-}$.
 Thus the optical solution is  Eq.(\ref{SolD}).
 Thus we have proved the claim.
 $\square$

\begin{theorem}
Let us assume that $\sigma\ne 0$ and $\delta \ne 0$.

Case 1 : Let the root  $\lambda$ be  that of Eq.(\ref{eq20a}). Then

1-1. if (1) $\delta > 0$ and $k>0$ or (2) $\delta < 0$ and $\eta_i>\eta_o$, the optical solution is Eq.(\ref{SolB}),

1-2. and if (1) $\delta > 0$ and $k<0$ or (2) $\delta < 0$ and $\eta_i < \eta_o$, the optical solution is Eq.(\ref{SolD}).

Case 2 : Let the root  $\lambda$ be  that of Eq.(\ref{eq20b}). Then

2-1. if (1) $\delta > 0$, $\sigma > 0$ and $\eta_i>\eta_o$ or (2) $\delta > 0$, $\sigma < 0$ and $k>0$ or (3) $\delta < 0$ and $k>0$, the optical solution is Eq.(\ref{SolB}),

2-2. and  if (1) $\delta > 0$, $\sigma > 0$ and $\eta_i < \eta_o$ or (2) $\delta > 0$, $\sigma < 0$ and $k<0$ or (3) $\delta < 0$ and $k<0$, the optical solution is Eq.(\ref{SolD}).

Case 3 : Let the root  $\lambda$ be  that of Eq.(\ref{eq20c}). Then

3-1. if (1) $\delta > 0$, $\sigma > 0$ and $k>0$ or (2) $\delta > 0$, $\sigma < 0$ and $\eta_i>\eta_o$ or (3) $\delta < 0$ and $k>0$, the optical solution is Eq.(\ref{SolB}),

3-2. and if (1) $\delta > 0$, $\sigma > 0$ and $k<0$ or (2) $\delta > 0$, $\sigma < 0$ and $\eta_i<\eta_o$ or (3) $\delta < 0$ and $k<0$, the optical solution is Eq.(\ref{SolD}).

\end{theorem}

{\bf Proof:} Before we prove the claim, let us observe that the  arccosine function returns $\theta$ with $0\le
\theta \le \pi$ in computation. Thus we have
\begin{eqnarray}
\cos(\frac{\theta+2\pi}{3})\le \cos(\frac{\theta-2\pi}{3}) \le
\cos(\frac{\theta}{3}),\label{inequal}
\end{eqnarray}
where the 1st equality is valid when $\theta=0$ and the 2nd equality
is when $\theta=\pi$.

 To prove our claim, it is enough to show which $\lambda_0$ in Eq.(\ref{eq27}) corresponds to the chosen $\lambda$ in Eq.(\ref{eq20}).
Hence   Eq.(\ref{inequal}) shows the order of the three roots $\lambda$ in
Eq.(\ref{eq20}). Now this order can be compared with that of
$\lambda_{0}$ in Eq.(\ref{eq27}).

It is trivial to see that $\lambda_{01} < \lambda_{02}$ since $k \ne 0$. On the other
hand, the difference between $\lambda_{02}$ and $\lambda_{03}$  is
\begin{eqnarray}
\lambda_{02}-\lambda_{03}&=&2k^2-\frac{2\epsilon^2m^2\eta_i(\eta_i-2\eta_o+m\eta_i)}{1+m} \\
&=&\frac{2\epsilon^2}{1+m}[\eta_o^2-2m\eta_i\eta_o+m\eta_o^2]=\frac{2\epsilon^2}{1+m}\delta
.\nonumber
\end{eqnarray}
Since $m>0$, we have
\begin{eqnarray}
&&\lambda_{01} < \lambda_{03} < \lambda_{02} \,\,\, \mbox{if}\,\,\, \delta>0 \,\,\mbox{and}\,\, \sigma>0,\\
&&\lambda_{03} < \lambda_{01} < \lambda_{02} \,\,\, \mbox{if}\,\,\, \delta>0 \,\,\mbox{and}\,\, \sigma<0,\nonumber\\
&&\lambda_{01} < \lambda_{02} < \lambda_{03} \,\,\, \mbox{if}\,\,\, \delta<0.\nonumber
\end{eqnarray}
Case 1: Let $\lambda$ be that of Eq.(\ref{eq20a}). Then $\lambda$ is the largest root of Eq.(\ref{lambda}) at $y=0$.
Hence $\lambda=\lambda_{02}$ if $\delta>0$, and  $\lambda=\lambda_{03}$ if $\delta<0$. Case 2: Let $\lambda$ be that of Eq.(\ref{eq20b}).
Then $\lambda=\lambda_{03}$ if $\delta>0$ and $\sigma>0$, and $\lambda=\lambda_{01}$ if $\delta>0$ and $\sigma<0$, and $\lambda=\lambda_{02}$ if $\delta<0$. Case 3: Let $\lambda$ be that of Eq.(\ref{eq20c}).
Then $\lambda$ is the smallest root of Eq.(\ref{lambda}) at $y=0$.
Thus $\lambda=\lambda_{01}$ if (1) $\delta>0$ and $\sigma>0$ or (2) $\delta<0$, and $\lambda=\lambda_{03}$ if $\delta>0$ and $\sigma<0$. Finally,  the claim is proved by Lemma 2. $\square$

\vspace{0.5cm}
{\bf Remark:}
 In the case when $\sigma=0$ or $\delta=0$, the criteria for the choice of the optical solution
 between Eq.(\ref{SolB}) and Eq.(\ref{SolD}) is as follows.
In fact, if $\sigma=0$,  $\lambda_{01}=\lambda_{03}$, which is
a double real root, and $\lambda_{02}$ is the largest root of Eq.(\ref{lambda}).
On the other hand, if $\delta=0$,  $\lambda_{02}=\lambda_{03}$,
a double real root, and $\lambda_{01}$ is the smallest root of Eq.(\ref{lambda}).
For the optical solution to be continuous through $\triangle = 0$,
we should choose the single root $\lambda_{02}$ or $\lambda_{01}$ for each case
since the single root is identical to the root in Eq.(\ref{eq16}) in the limits.
As seen in Eq.(\ref{eq28}), $\lambda_{22}$ and $\lambda_{21}$ corresponding to
$\lambda_{02}$ and $\lambda_{01}$ respectively are still valid
only for each $\sigma=0$ and $\delta=0$ case.
Thus the proof of Lemma 2 is still valid for each case. Hence the optical solution is Eq.(\ref{SolB}) if $k>0$.
It is Eq.(\ref{SolD}) if $k<0$. In practical computation, however,
it would be better to use the root $\lambda$ in Eq.(\ref{eq16}) initially
 when $\sigma=0$ or $\delta=0$.

In the next subsection, we shall discuss more on the continuities of optical solutions.

\subsection{ The continuities and interpolations of the optical solutions}

Let us suppose that $\delta\ne0$ and $\sigma \ne 0$.
We note that $\triangle$ in Eq.(\ref{triangle}) is a cubic polynomial in $y^2$.
Thus we may apply the previous cubic equation solving method to the equation $\triangle=0$ again.
Furthermore, what makes things simple is that
not only $\triangle<0$ when $y=0$ but also $D_6<0$ in Eq.(\ref{triangle_coeff}).
Thus either there is no positive solution $y$
or there is a pair of positive solutions such that $\triangle=0$.

In fact, from $\tilde q \equiv (3 D_6 D_2-D_4^2)/(3 D_6^2)$,
$\tilde r \equiv (9 D_6 D_4 D_2-27 D_6^2D_0-2D_4^3)/(27 D_6^3)$ and $\tilde \triangle \equiv \tilde r^2+\frac{4}{27}\tilde q^3$
corresponding to Eq.(\ref{eqdelta}), it follows that $\tilde \triangle > 0$ means that there is no solution $y$
such that $\triangle=0$. However, if $\tilde \triangle < 0$,
either there is no positive solution $y$ such that $\triangle=0$,
or there are only two positive solutions $y$.

If there is no solution $y$ such that $\triangle=0$, then any initially chosen root $\lambda$ makes the optical solution
continuous for all $y$ on its domain.
However, if there is a positive solution $y$,
the initially chosen root $\lambda$ should be changed as $y$ increases for the continuity of the optical solution.
Now we want to choose the root $\lambda$ such that it is continuous through $\triangle=0$ when
$\triangle$ changes its sign.
The following claim is very useful in practical computation.
\begin{lemma}
Let $\sigma\ne 0$, $\delta \ne 0$ and
 $y_0$ be a positive solution such that $\triangle=0$.

1. If $r>0$ at $y_0$, the only continuous root through $\triangle=0$ is  Eq.(\ref{eq20a}).

2. If $r<0$ at $y_0$, the only continuous root through $\triangle=0$ is  Eq.(\ref{eq20c}).
\end{lemma}
{\bf Proof:} For notational convenience, let us write $\lambda^{(0)}$, $\lambda^{(1)}$, $\lambda^{(2)}$ and $\lambda^{(3)}$
respectively for the root $\lambda$'s of Eq.(\ref{eq16}), Eq.(\ref{eq20a}),  Eq.(\ref{eq20b}) and  Eq.(\ref{eq20c}).
If $r>0$ at $y_0$, $\theta=0$ or $r=2\rho$.
Now it is easy to see that $\lim_{\triangle \to 0^-}\lambda^{(1)}= -b/(3a)+2\rho^{\frac{1}{3}}= \lim_{\triangle \to 0^+}\lambda^{(0)}$,
while $\lim_{\triangle \to 0^-}\lambda^{(2)}=\lim_{\triangle \to 0^-}\lambda^{(3)}=-b/(3a)-\rho^{\frac{1}{3}}$.
On the other hand,
If $r<0$ at $y_0$, $\theta=\pi$ or $r=-2\rho$. Thus
$\lim_{\triangle \to 0^-}\lambda^{(3)}= -b/(3a)-2\rho^{\frac{1}{3}}= \lim_{\triangle \to 0^+}\lambda^{(0)}$,
while $\lim_{\triangle \to 0^-}\lambda^{(1)}=\lim_{\triangle \to 0^-}\lambda^{(2)}=-b/(3a)+\rho^{\frac{1}{3}}$.
 $\square$

\vspace{0.5cm}
{\bf Remark:} In the case when there is a positive solution $y$ such that $\triangle=0$,
we may use the criteria given by Lemma 3 and Theorem 1 together
for the initial choices of the root $\lambda$ and the optical solution
of a Cartesian oval. As $y$ increases,
the initially chosen root $\lambda$ should be changed for the continuity of the optical solution.
In fact, there is no choice except the root $\lambda$ of Eq.(\ref{eq16}) for the transition
 from $\triangle <0$ to $\triangle >0$.
On the other hand, if the transition  direction is reversed, i.e. from $\triangle >0$ to $\triangle <0$,
there can be three choices in Eq.(\ref{eq20}) for  the root $\lambda$.
However, we should choose the root $\lambda$ according to the above Lemma 3
 in order to make the optical solution continuous.

\vspace{0.5cm}

Furthermore, for a continuous optical solution in Eq.(\ref{SolB}) or Eq.(\ref{SolD})
\begin{eqnarray}
z=B/[1+(1-AB)^{\frac{1}{2}}]\, \,\, \mbox{  }\,\, (\, A\equiv A_{\pm} \,\,\mbox{and} \,\, B\equiv B_{\pm}\,), \label{TheOptSol}
\end{eqnarray}
it might happen that $1-AB<0$ as $y$ increases.
In this case, we can interpolate it continuously by putting $1-AB=0$ in a manner similar to conic case in optical design.
That is, we may interpolate the optical solution in  Eq.(\ref{TheOptSol}) by the curve $z=B$ in the region where $1-AB<0$.
\begin{definition}
Let $z=B/[1+(1-AB)^{\frac{1}{2}}]$ be an optical solution. Then
$z=B$ is called an interpolating curve of the optical solution  in the region where $1-AB<0$.
\end{definition}

We conclude this section with the observation that the coefficient
of the second order term in $y$ of the optical solution is  the
curvature $c_0$ of the optical solution from Eqs.(\ref{approxiB}, \ref{approxiD})
\begin{eqnarray}
c_0&\equiv&\frac{ \eta_i-m\eta_o}{ \epsilon (1-m)}.\label{curvature}
\end{eqnarray}
When $y$ is small, it is also  interesting to see that  the optical
solution is of the following form.
\begin{eqnarray}
z&=&\frac{c_0y^2+O(y^4)}{1+\sqrt{1-(1+K)c_0^2y^2+O(y^4)}}
\end{eqnarray}
for some constant $K$ as a function of constants $\eta_i, \eta_o$ and $m$. This form for the optical solution of a
Cartesian oval shows its deviation  from a conic curve when $y$ is
small. It looks like that $K$ plays the role of a conic constant.
However, the insightful relations of the optical solutions to conics
shall be shown in the next section.

\section{Superconics as extensions of aspheric curves based on conics }

It is well-known in the literature\cite{Bev2006, Law1972, Loc1967} that  Cartesian ovals become conics if $m=\pm 1$,
which can be easily observed in Eq.(\ref{eq2}) of this work.
In this section, however, we shall consider the limits of the optical solution expressed
in Eq.(\ref{SolB}) or Eq.(\ref{SolD}) from a different point of view.

In fact, we want to see the limits of both  optical solutions and their interpolating curves.
First of all,  we observe that Eq.(\ref{eq2}) is invariant under the
replacements of  $\eta_i$, $\eta_o$ and $m$ by $\eta_o$, $\eta_i$
and $1/m$ respectively:
\begin{eqnarray}
\eta_i &\Rightarrow& \eta_o, \,\,\, \eta_o \Rightarrow \eta_i, \,\,\,
m \Rightarrow   1/m. \label{replacement}
\end{eqnarray}
Thus any optical solution with $0<m<1$ may be represented by an
optical solution with $m>1$.  From now on, we may assume
that $m>1$ without loss of generality.

\begin{theorem}
Let $m>1 $ and $z=B/[1+(1-AB)^{\frac{1}{2}}]$ be an optical solution. Then we have
\begin{align}
\lim_{\eta_i \rightarrow 0} \frac{B}{1+(1-AB)^{\frac{1}{2}}} &= \frac{c_0y^2}{1+[1-(1+K)c_0^2y^2]^{\frac{1}{2}}}
&\mbox{if} \,\,\,1-AB \ge 0 ,\\
\lim_{\eta_i \rightarrow 0} B &= c_0y^2
&\mbox{if}\,\,\, 1-AB < 0 , \nonumber
\end{align}
where the curvature $c_0$ and the conic constant $K$ are defined as
\begin{eqnarray}
c_0&\equiv&\frac{\eta_o}{\epsilon(1-\frac{1}{m})}, \hspace{1cm}
K\equiv -\frac{1}{m^2}.\label{eq52}
\end{eqnarray}
\end{theorem}

{\bf Proof:} Let us choose the root $\lambda$ of Eq.(\ref{eq20a}) to make the optical solution continuous
for all $y$ as can be seen in the below.  Now for small $\eta_i$ we put
\begin{eqnarray}
\rho^\frac{1}{3}&=&\sqrt{-\frac{q}{3}}\equiv
[\rho^\frac{1}{3},0]+[\rho^\frac{1}{3},1]\eta_i +
[\rho^\frac{1}{3},2]\eta_i^2,\nonumber
\end{eqnarray}
where the notation $[f, n]$ represents the $n$th order coefficient of $f$ in $\eta_i$.
Then it is not hard to find
\begin{eqnarray}
\,[\rho^\frac{1}{3},0]&=& \frac{2}{3}\epsilon^2\eta_o^2,\label{rhoco}\\
\,[\rho^\frac{1}{3},1]&=& -\frac{2(2+m)m\epsilon^2\eta_o}{3(1+m)},\nonumber\\
\,[\rho^\frac{1}{3},2]&=&\frac{\epsilon^2m^2(1+2m+4m^2)}{3(1+m)^2}-\frac{2}{3}(1+2m^2)\eta_o^2y^2
.\nonumber
\end{eqnarray}

Now we put
\begin{eqnarray}
\cos\frac{\theta}{3}&\equiv&[\cos \frac{\theta}{3},0]+[\cos
\frac{\theta}{3},1]\eta_i+[\cos \frac{\theta}{3},2]\eta_i^2.
\end{eqnarray}
From the observations that $\cos \theta = 4\cos^3\frac{\theta}{3}-3\cos\frac{\theta}{3}$  and
\begin{eqnarray}
\,[\cos\theta,0]&=& 1,\\
\,[\cos\theta,1]&=& 0,\nonumber\\
\,[\cos\theta,2]&=& -\frac{27\epsilon^4\eta_o^4m^4[\epsilon^2(m-1)-(1+m)\eta_o^2y^2]}{2(m-1)(1+m)^2\epsilon^6\eta_o^6},\nonumber
\end{eqnarray}
we have
\begin{eqnarray}
\,[\cos \frac{\theta}{3},0]&=& 1,\label{thetaco}\\
\,[\cos \frac{\theta}{3},1]&=& 0, \nonumber\\
\,[\cos \frac{\theta}{3},2]&=&  -\frac{3m^4}{2(1+m)^2\eta_o^2}
+\frac{3m^4y^2}{2\epsilon^2(m-1)(1+m)}.\nonumber
\end{eqnarray}
Thus if we put
\begin{eqnarray}
\lambda
&\equiv&[\lambda,0]+[\lambda,1]\eta_i+[\lambda,2]\eta_i^2,\nonumber
\end{eqnarray}
 it is straightforward from Eq.(\ref{eq20a}) and Eqs.(\ref{rhoco}, \ref{thetaco}) to see the followings
\begin{eqnarray}
\,[\lambda,0]&=&2\epsilon^2\eta_o^2 \label{lamco}\\
\,[\lambda,1]&=&-4m\epsilon^2\eta_o.\nonumber\\
\,[\lambda,2]&=&2\epsilon^2m^2+ \frac{2}{m^2-1}\eta_o^2y^2.\nonumber
\end{eqnarray}
For $y=0$,   $ \lambda \approx
2\epsilon^2\eta_o^2-4m\epsilon^2\eta_o\eta_i+2\epsilon^2m^2\eta_i^2
=2k^2$ which corresponds to $\lambda_{02}$ in Eq.(\ref{lambda02}). Now we
put $A = A_0+A_1\eta_i+\cdots $, and $B = B_0+B_1\eta_i+\cdots $ and
suppose that $\eta_o
>0$. Then since $\delta\approx(1+m)\eta_o^2>0$ and $k=\epsilon(\eta_o-m\eta_i)=\epsilon\eta_o>0$,
 the use of  $A=A_{+}$ and $B=B_{+}$ according to the criteria in Theorem 1 yields
\begin{eqnarray}
A_0&=&\frac{(1+m)\eta_o}{\epsilon m }=(1+K)c_0,\,\,\,\, B_0
=\frac{m\eta_oy^2}{\epsilon(m-1)}=c_0y^2. \label{eq98}
\end{eqnarray}
On the other hand,   suppose that $\eta_o <0$. Then $k<0$ and we
obtain the same result for $A=A_{-}$ and $B=B_{-}$.
Hence the case when $1-AB \ge 0$ has been proved.

Now we observe  from Eq.(\ref{triangle}) that, when $\eta_i \rightarrow 0$,
\begin{eqnarray}
\triangle &\approx & -\frac{256\epsilon^{10}m^4\eta_o^{10}}{27(1+m)^2(m-1)}((m-1)\epsilon^2-(1+m)\eta_o^2y^2)\eta_i^2.
\end{eqnarray}
Thus the root of Eq.(\ref{eq20}) for $\triangle<0$ is changed once to that of Eq.(\ref{eq16}) at $y=\sqrt{\frac{(m-1)\epsilon^2}{(m+1)\eta_o^2}}=\sqrt{\frac{1}{(1+K)c_0^2}}$
through $\triangle=0$.
Hence it remains to show the limit when $\triangle>0$.

When $\eta_i \rightarrow 0$, let us put $\lambda \approx
\lambda_0+\lambda_1\eta_i+\lambda_2\eta_i^2 $ again for $\triangle>0$. Then from
Eq.(\ref{eq16}), the same coefficients as in Eq.(\ref{lamco}) are obtained by a little bit
lengthy calculation. Hence we have $B \rightarrow c_0y^2$
for $y\ge \sqrt{\frac{1}{(1+K)c_0^2}}$ or $1-AB=1-(1+K)c_0^2y^2 \le 0$.
  $\Box$

\vspace{0.5cm}

 Any
central conic with $c_0, K$($0<K<1$) determines $\eta_o, m$ by Eq.(\ref{eq52}).
We remark that the circle is obtained in the limit when $m\rightarrow
\infty$ and the parabola is the limit of the central conic when
$m\rightarrow 1$ with $c_0$ being fixed.

\begin{theorem}
Let $m>1$ and $z=B/[1+(1-AB)^{\frac{1}{2}}]$ be an optical solution. Then we have
\begin{eqnarray}
&&\lim_{\eta_o \rightarrow 0} \frac{B}{1+(1-AB)^{\frac{1}{2}}} = \frac{c_0y^2}{1+[1-(1+K)c_0^2y^2]^{\frac{1}{2}}}
\,\,\,\mbox{for all} \,\,\,y,   
\end{eqnarray}
where the curvature $c_0$ and the conic constant $K$ are defined as
\begin{eqnarray}
c_0&\equiv&\frac{\eta_i}{\epsilon(1-m)}, \hspace{1cm}
K\equiv -m^2.\label{eq54}
\end{eqnarray}
\end{theorem}

{\bf Proof:} On the contrary to the case when $\eta_i \rightarrow 0$,
\begin{eqnarray}
\triangle &\approx & -\frac{256\epsilon^{10}m^{10}\eta_i^{10}}{27(1+m)^2(m-1)}((m-1)\epsilon^2+(1+m)\eta_i^2y^2)\eta_o^2,
\end{eqnarray}
when  $\eta_o \rightarrow 0$. Thus $\triangle<0$ for all $y$ since $m>1$.
That is, there is no solution $y$ such that $\triangle=0$.
Hence we may choose the root of Eq.(\ref{eq20a}) for the convenience of computation.

In a similar manner as in the previous theorem,
 we put
\begin{eqnarray}
\rho^\frac{1}{3}&=&\sqrt{-\frac{q}{3}}\equiv
[\rho^\frac{1}{3},0]+[\rho^\frac{1}{3},1]\eta_o +
[\rho^\frac{1}{3},2]\eta_o^2,\nonumber
\end{eqnarray}
where the notation $[f, n]$ represents the $n$th order coefficient of $f$ in $\eta_o$.
Then it is not hard to find
\begin{eqnarray}
\,[\rho^\frac{1}{3},0]&=& \frac{2}{3}\epsilon^2m^2\eta_i^2, \label{rhoco1}\\
\,[\rho^\frac{1}{3},1]&=&
-\frac{2(1+2m)m\epsilon^2\eta_i}{3(1+m)},\nonumber\\
\,[\rho^\frac{1}{3},2]&=&\frac{\epsilon^2(4+2m+m^2)}{3(1+m)^2}-\frac{2}{3}(2+m^2)\eta_i^2y^2
.\nonumber
\end{eqnarray}

Now we observe that
\begin{eqnarray}
\,[\cos\theta,0]&=& -1,\\
\,[\cos\theta,1]&=& 0,\nonumber\\
\,[\cos\theta,2]&=& Q,\nonumber
\end{eqnarray}
where for $m>1$
\begin{eqnarray}
Q&\equiv& \frac{27}{2m^2(1+m)^2\eta_i^2}
+\frac{27y^2}{2\epsilon^2m^2(m^2-1)}>0.\label{eqQ}
\end{eqnarray}
We put
\begin{eqnarray}
\cos\frac{\theta}{3}&\equiv&[\cos \frac{\theta}{3},0]+[\cos
\frac{\theta}{3},1]\eta_o+[\cos \frac{\theta}{3},2]\eta_o^2.
\end{eqnarray}
Then from the observations that $\cos \theta = 4\cos^3\frac{\theta}{3}-3\cos\frac{\theta}{3}$  and
$\cos (\arccos(-1+cx^2)/3)= \frac{1}{2}+\sqrt{\frac{c}{6}}x
-\frac{c}{18}x^2+\cdots$ for $ c\ge 0$ and some $x$,
it follows that
\begin{eqnarray}
\,[\cos \frac{\theta}{3},0]&=& \frac{1}{2},\label{thetaco1}\\
\,[\cos \frac{\theta}{3},1]&=& \sqrt{\frac{Q}{6}},\nonumber\\
\,[\cos \frac{\theta}{3},2]&=& -\frac{Q}{18}.\nonumber
\end{eqnarray}

Thus if we put
\begin{eqnarray}
\lambda
&\equiv&[\lambda,0]+[\lambda,1]\eta_o+[\lambda,2]\eta_o^2,\nonumber
\end{eqnarray}
 it is straightforward from Eq.(\ref{eq20a}) and Eqs.(\ref{rhoco1}, \ref{thetaco1}) to see the followings
\begin{eqnarray}
\,[\lambda, 0]&=&2\epsilon^2m^2\eta_i^2,\,\,\,\,
[\lambda,1]=-\frac{2(1+2m)m\epsilon^2\eta_i}{1+m} +\frac{
4\epsilon^2m^2\eta_i^2}{3\sqrt{6}}\sqrt{Q},
\\
\,[\lambda,2]&=&-\frac{2(m^2-1)(\epsilon^2+(m^2-1)\eta_i^2y^2)}{3(1-m^2)}
+\epsilon^2m^2\eta_i^2\cdot[\frac{1}{m^2(1+m)^2\eta_i^2}
+\frac{y^2}{\epsilon^2m^2(m^2-1)}]\nonumber\\
&&+[\frac{\epsilon^2(4+2m+m^2)}{3(1+m)^2}-\frac{2}{3}(2+m^2)\eta_i^2y^2]
-\frac{2(1+2m)m\epsilon^2\eta_i}{3(1+m)}\cdot \frac{
1}{\sqrt{6}}\sqrt{Q}.\nonumber
\end{eqnarray}

The optical solution is Eq.(\ref{SolB}) as can be seen in the below. We
 expand $A_{+}$ and $B_{+}$ in terms of $\eta_o$: $A = A_0+A_1\eta_o+\cdots $,
and $B = B_0+B_1\eta_o+\cdots $.  Then it is  lengthy but straightforward  to obtain the following.
\begin{eqnarray}
A_{0}&=&0,\,\,\, A_{1} = \frac{m\pm 1}{\epsilon m },\,\,\,
B_0=\frac{2\epsilon}{(1+m)\eta_i}\pm\frac{2\sqrt{6}\epsilon
m\sqrt{Q}}{9},
\end{eqnarray}
where $(+)$ and $(-)$ are for $\eta_i<0$  and $\eta_i > 0$ respectively.
Thus  we have
\begin{eqnarray}
\lim_{\eta_o \rightarrow 0} \frac{B}{1+(1-AB)^{\frac{1}{2}}}&=& \frac{B_0}{2} =\frac{ \epsilon}{ (1+m)\eta_i}\pm\frac{\sqrt{6}\epsilon m\sqrt{Q}}{9}\label{eq61}\\
&=& \frac{1}{(1+K)c_0}\pm \frac{\epsilon}{(1+m)|\eta_i|}\sqrt{1-(1+K)c_0^2y^2}\nonumber\\
&=& \frac{c_0y^2}{1+[1-(1+K)c_0^2y^2]^{\frac{1}{2}}}.\nonumber
\end{eqnarray}
Thus the claim has been proved.  $\square$

Any central conic with $c_0, K$($K>1$) determines $\eta_i, m$ by Eq.(\ref{eq54}). In this
case, the parabola is the limit of the central conic when $m\rightarrow
1$ with $c_0$ being fixed.

\begin{definition}
 A superconic curve is defined to be an aspheric curve based on
optical solution described as follows: For any $y$,
\begin{eqnarray}
z=B/[1+(1-AB)^{\frac{1}{2}}]+ \sum_{n=2}^{N} f_{2n} y^{2n},\label{asp2}
\end{eqnarray}
where $N$ is a positive integer and $f_{2n}$'s are some constants.
\end{definition}
The superconic curve  is obviously different from Greynolds'\cite{Ala2002}.
Now  we have the main result of this work, which follows obviously from the previous theorems.

\begin{corollary}
A superconic curve described by Eq.(\ref{asp2})
is an extension of an aspheric curve based on conic described by Eq.(\ref{asp}).
\end{corollary}

\section{Example : A Family of Optical Solutions and Conics}

We note that $m$($m>1$) and $\eta_i, \eta_o$ are free variables for the
optical solutions of  Cartesian ovals. However, there may be some
other choices for free variables. It seems to be nice to choose the curvature $c_0$,
$m$ and $\eta_i$ as free variables, in which case $\eta_o$ can be obtained by
$\eta_o = (\eta_i-\epsilon (1-m)c_0)/m$ from Eq.(\ref{curvature}).

With this choice of free variables, we may describe the set of
 curves with the same curvature $c_0$ and $m$ but different $\eta_i$'s as one family.
Of course, typical members of the family are the optical solutions with such $c_0$ and $m$.
However, it is interesting  that there are some special members in the family.
 That is, if $\eta_i=0$ or
$\eta_i=\epsilon(1-m)c_0$ (i.e. $\eta_o=0$), the curves become
conics. In fact, when $\eta_i=0$, the curve is an ellipse with the curvature $c_0$ and the conic constant
$K=-1/m^2$. When  $\eta_i=\epsilon(1-m)c_0$, the curve is a hyperbola with the same curvature $c_0$ and the conic constant
$K=-m^2$.

Especially,   as can be seen from Eq.(\ref{eq2}), if $\eta_i=\epsilon c_0$ (i.e. $\eta_i=\eta_o$)
or $\eta_i=\epsilon c_0/(1+m)$ (i.e. $\eta_i=\eta_o/m$ or $k=0$), they represent one
circle with the same curvature $c_0$ and the conic constant $K=0$.
It may be easily observed that it  is the circle  possessed
 commonly  by all the families with the
same curvature $c_0$.
Thus even conic curves with such special $\eta_i$'s in the above
 may be regarded as members of the family of curves with $c_o$ and $m$.

In Fig. 1, we demonstrate an example of a
family of curves with the same
curvature $c_0=0.3$ and $m=1.5$ for  $\epsilon=0.6$.
The family contains the optical solutions with the same $c_0$ and $m$
 but with different $\eta_i$'s in
the range $0.3(C_1) \sim -0.3(C_6)$ following the arrow.

Moreover, the family contains the conics as well with the same
curvature $c_0$ and $m$ but with different $\eta_i$'s in
the range $0.18(D_1) \sim -0.09(D_4)$. Here we note that the conics $D_1$ and $D_2$ are the same circle
with  $\eta_i=\epsilon c_0=0.18$
and $\eta_i=\epsilon c_0/(1+m)=0.072$  respectively, which is a peculiar fact. It is observed that
the optical solutions, e.g. $C_2$,  with $\eta_i$'s between
$0.18$ and $0.072$ are very close to the circle.

\begin{figure} 
\includegraphics[scale=1.0]{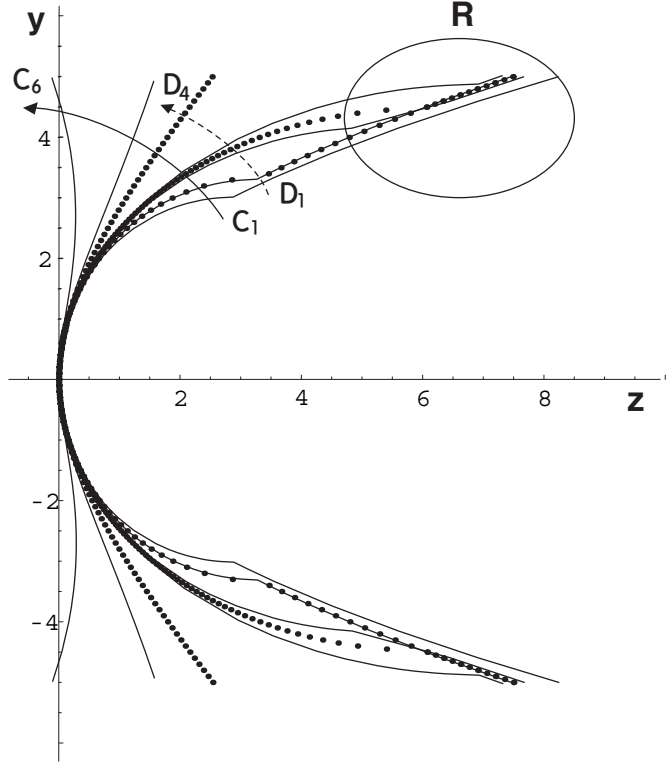}
\caption{The Figure shows a family of curves containing of optical solutions
$C_1, \cdots, C_6$ (real lines) and conics $D_1(=D_2)$, $D_3$ and $D_4$ (dotted lines) following the arrows respectively with the same curvature
$c_0=0.3$ and $m=1.5$, but with different $\eta_i$'s.
The values of $\eta_i$'s of the curves are $\eta_i=0.3 (C_1)$, $0.18(D_1)$, $0.15(C_2)$, $0.072(D_2)$, $0.01(C_3)$,
$0(D_3)$, $-0.01(C_4)$, $-0.09(D_4)$, $-0.15(C_5)$ and $-0.3(C_6)$ respectively.
The interpolating curves of the conics $D_1, D_2$ and $D_3$ are  $z=c_0y^2$ when $1-(1+K)c_0^2y^2 <0$.
Similarly, the interpolating curves of the optical solutions $C_1, \cdots, C_4$ are $z=B$ when $1-AB<0$. They are drawn in the region $R$.
}
\end{figure}

The curve $D_3$ is an ellipse with $\eta_i=0$. The conic constant is $K=-1/m^2 \approx -0.44$.
The optical solutions with
small $\eta_i$ such as $C_3$ or $C_4$ are slightly deviated from the
conic curve $D_3$. If $\eta_i$ is not small, e.g. $C_1, C_6$, it is
far different from the conic curve although it has the same
 $c_0$ and $m$.

 The curve $D_4$ is
the case when $\eta_i=\epsilon(1-m)c_0=-0.09$. It is
interesting to note that $D_4$ is a hyperbola with the
 conic constant $K=-m^2=-2.25$ by Eq.(\ref{eq54}) since $\eta_o= 0$. In fact,
there always exists such a pair of conics as an ellipse $D_3$ and a
hyperbola $D_4$ for the family with a given $c_0$ and $m$, and both of them approach one
parabola when $K\rightarrow -1$ (i.e. $m \rightarrow 1$).

The interpolating curves of the optical solutions are described by $z=B$
like the interpolating conic curves are described by $z=c_0y^2$ for the case
when $1-(1+K)c_0^2y^2 <0$  in Eq.(\ref{conicEx}).
The interpolating curves are shown in the region $R$.

For most of the optical solutions in Fig.1, the initially chosen root $\lambda$ of Eq.(\ref{eq20a}) has been changed
to that of Eq.(\ref{eq16}) at some $y$. Furthermore, for the curves $C_1$ and $C_2$, the root $\lambda$
has been changed again even to that of
Eq.(\ref{eq20c}) in the region $R$.

It is natural to extend the family of optical solutions and conics to the family of superconic curves
and aspheric curves based on conics by adding the higher order polynomial terms $\sum_{n=2}^{N} f_{2n} y^{2n}$
to optical solutions and conics as in Eq.(\ref{asp2}) and  Eq.(\ref{asp}).

The optical solutions may be used as the
starting curves for lenses at the initial design of an optical system using the perfectly
focusing property.  And the initial optical design may be followed by an elaborate optimization, in which case
 superconic curves in Eq.(\ref{asp2}) seem to be the most suitable curves for the optimization since they are extensions of not only the initial optical solutions but also aspheric curves based on conics in Eq.(\ref{asp}).

\section{CONCLUSIONS}

In this work, we have investigated the criteria to find the optical solution of a Cartesian oval
and discussed the continuity and interpolation of the optical solution. Moreover, we have shown that conics and their interpolating curves are the limiting cases of the optical solutions and their interpolating curves respectively. It follows then that all those curves with the same curvature $c_0$ and parameter $m$ but different $\eta_i$'s, including both optical solutions and conics,  can be regarded as members of one family of curves. We have demonstrated
an example about a family of curves.

Most of all, the above work on the optical solutions makes it possible to construct
another kind of superconic curves in Eq.(\ref{asp2}) that are different from Greynolds'. That is,
the superconic curves suggested in this work are  extensions of not only  conics and Cartesian ovals
 but also aspheric curves based on conics in Eq.(\ref{asp}) while Greynolds' superconic curves are not extensions
 of aspheric curves based on conics. Also the superconic curves in this work are expressed
 in explicit form as in Eq.(\ref{asp2}) while   Greynolds' superconic curves are in implicit form as in Eq.(\ref{Reynolds}).

The relationship between superconic curves in Eq.(\ref{asp2}) and aspheric curves based on conics in Eq.(\ref{asp}) seems
to be a promising property  that Greynolds'
superconic curves do not have.

\vspace{0.5cm}
{\bf Acknowledgments}
\vspace{0.5cm}

The author would like to thank Professor B. S. Lee  for helpful discussions and critical comments.
This paper was supported by the Semyung University Research Grant of 2013.

\end{document}